\documentclass[11pt]{article}
\usepackage{amsmath,amsbsy,amssymb,natbib,fullpage,
           latexsym,graphicx}

\newtheorem{prop}{Proposition}
\newcommand{\bm}[1]{\mbox{\boldmath{$#1$}}}
\newcommand{\mc}[1]{\ensuremath{{\mathcal #1}}}
\newcommand{\m}[1]{\mbox{\bf{#1}} }
\begin{document}

\title{ 
Model averaging  and 
dimension selection for the singular value decomposition }
\author{ Peter D. Hoff  
\thanks{ Departments of Statistics, Biostatistics and the
Center for Statistics and the Social Sciences,
University of Washington,
Seattle, Washington 98195-4322, U.S.A..
Email: {\tt hoff@stat.washington.edu}.
Web: {\tt http://www.stat.washington.edu/hoff/}. 
This research was supported by Office of Naval Research grant
N00014-02-1-1011 and National Science Foundation grant
SES-0417559. The author thanks David Hoff 
 and Michael Perlman for helpful discussions.
        } 
       }
\date{ \today }
\maketitle
    
\begin{abstract}
Many multivariate data analysis techniques for  
an $m\times n$ 
matrix $\m Y$ are related to the 
model $\m Y = \m M +\m E$, where $\m Y$ is an $m\times n$ matrix of 
full rank and $\m M$ is an unobserved mean matrix of rank $K< (m\wedge n)$. 
Typically the rank of $\m M$ is estimated in a heuristic way and then the 
least-squares 
estimate of $\m M$ is obtained via the singular value decomposition of 
$\m Y$, yielding an estimate that can have a very high variance. 
In this paper we suggest a model-based alternative to the above approach
by providing prior distributions and posterior estimation for
the rank of $\m M$ and the components of its singular 
value decomposition. 
In addition to providing more accurate inference, such an approach 
has the advantage of being extendable to more general data-analysis 
situations, such as inference in the presence of missing data and 
estimation in a generalized linear modeling framework. 

\vspace{.2in}
\noindent {\it Some key words}: 
interaction, 
model selection, multiplicative effects, multiple hypergeometric function,
relational data, social network, Steifel manifold. 
\end{abstract}

\section{Introduction}
Many inferential and descriptive methods for
multivariate and matrix-valued data 
are variations on the idea of modeling the 
data matrix $\m Y$ as equal to 
 a reduced-rank mean matrix $\m M$ plus Gaussian
noise, and  estimating $\m M$  after deciding upon its rank. 
A concept that is central to such models is the singular value 
decomposition: 
every $m\times n$ matrix $\m M$ has a representation of the form 
$\m M= \m U\m D \m V'$
where, in the case $m\geq n$, 
$\m U$ is an $m\times n$ matrix with orthonormal columns, 
$\m V$ is an $n\times n$ matrix with orthonormal columns and
$\m D$ is an $n\times n$  diagonal matrix, with diagonal elements
     $\{ d_1,\ldots, d_n\}$ typically taken to be a decreasing sequence  of 
    non-negative numbers. 
The triple $\{\m U,\m D, \m V\}$ is called the singular value decomposition of $\m M$. 
The squared elements of the diagonal of $\m D$ are the eigenvalues of
$\m M'\m M$ and the columns of
$\m V$  are the corresponding eigenvectors.
The matrix $\m U$ can be obtained from the first $n$ eigenvectors of
$\m M \m M'$. The number of non-zero elements of $\m D$ is the rank of $\m M$. 

The appeal of the singular value decomposition is partly due
to its 
interpretation as a multiplicative model based 
on row and column factors. Given a model of the 
form $\m Y = \m M + \m E$, 
the  elements of $\m Y$ can be written
$y_{i,j} = \m u_i '\m D \m v_j  + e_{i,j}$, 
where $\m u_i$ and $\m v_j$ are 
the $i$th and $j$th rows of $\m U$ and  $\m V$ respectively. 
This model thus provides a representation of the
 systematic variation among the entries of $\m Y$
by row and column factors. 
Models of this type play a role in
the analysis of relational data (Harshman et al., 1982), 
biplots (Gabriel 1971, Gower and Hand 1996) and in 
reduced-rank interaction models for factorial designs (Gabriel 1978, 1998). 
The model is also closely related to factor analysis, in which the row 
vectors of $\m Y$ are modeled as i.i.d.\ samples from the model 
$\m y_i = \m u_i \m D \m V' + \m e_i$. 
In this situation, the 
goal is typically to represent the covariance across the $n$ columns
 by
their relationship to $K<n$ unobserved latent factors. 
Lawley and Maxwell (1971) provide a comprehensive overview of 
factor analysis models, as well as methods of maximum likelihood 
estimation for model parameters.

The singular value decomposition also plays a role in parameter estimation
for the reduced-rank model:
Assuming the rank of the mean matrix $\m M$ is $K<n$ and 
letting $(\hat {\m U},\hat {\m D},\hat {\m V})$ be the singular value decomposition 
of the data matrix $\m Y$, 
the least-squares estimate of $\m M$ (and maximum likelihood estimate 
under Gaussian noise) is given by 
$\hat {\m M}_K =  \hat {\m U}_{[,1:K]} \hat{\m D}_{[1:K,1:K]} \hat{\m V}_{[,1:K]}'$, 
where $\hat {\m U}_{[,1:K]}$ denotes the first $K$ columns of $\hat {\m U}$ and 
$\hat{ \m D}_{[1:K,1:K]}$ denotes the first $K$ rows and columns of $\hat {\m D}$
(Householder and Young 1938, Gabriel 1978). 
In applications such as signal processing, image analysis, and more recently
large-scale gene expression data, representing a noisy data matrix $\m Y$ 
by $\hat {\m M}_K$ with  $K<<n$ has the effect of capturing the main patterns of 
$\m Y$
while eliminating much of the noise.

Despite its utility and simplicity, several issues limit the use of 
the singular value decomposition as an estimation procedure. 
The first is that the rank $K$ of the approximating mean
 matrix  $\hat {\m M}_K$ must be specified. 
Standard practice is to plot the singular values of $\m Y$ in 
decreasing order and then select $K$ to be the index where the last 
``large gap,'' or ``elbow'' occurs. 
The second issue is that, even if the rank is chosen correctly, 
one may be concerned with its variance: 
For high-dimensional parameter spaces, 
the mean squared error of least-squares estimates is often 
higher than that of penalized  or Bayes estimates. 
Finally, non-model-based approaches are somewhat limited 
in their ability to handle missing or non-normal data.

Philosophical debates aside,
Bayesian methods often provide sensible procedures 
for model selection and  high-dimensional parameter estimation. 
For the model described above, 
a Bayesian procedure would provide a mapping from a prior 
distribution $p(\m U,\m D,\m V,\sigma^2 )$ to a posterior 
distribution $p(\m U,\m D,\m V,\sigma^2 |\m Y)$. 
Of primary interest
might be functions of this posterior distribution, such as $E(\m M|\m Y)$ or 
the marginal posterior distribution of the rank $p(K|\m Y) \propto
p(K)p(\m Y|K)$.  Both of these quantities require  integration over the 
complicated, high-dimensional space of $\{\m U,\m D,\m V\}$ 
 for each value of $K$. 
In the related factor analysis model where
the  elements of $\m U$ are modeled as 
independent normal random variables, 
the difficulty in calculating marginal probabilities
 has led to the development of approximate Bayesian procedures:
Rajan and Rayner (1997) provide
a coarse approximation to the marginal probability $p(\m Y| K)$ by plugging 
in maximum-likelihood estimates. Minka (2000) improves on this 
by providing 
a Laplace approximation to the desired marginal probability.
Both of these procedures rely on asymptotic approximations, 
and do not provide 
Bayesian estimates of $\m M$ once the dimension has been selected. 

Alternatively one could turn to Markov chain Monte Carlo methods 
to obtain approximations to the integrals that are necessary
for model selection.  Tierney (1994) describes the theory of MCMC 
for very general parameter spaces, 
and Green (1995) outlines conditions under which the
reversibility of the Metropolis-Hastings algorithm is 
maintained for model-selection
problems. However, obtaining efficient 
proposal distributions for high-dimensional
 problems is a difficult and delicate art. 
For the factor analysis problem, 
Lopes and West (2004) devise a workable reversible-jump
algorithm by constructing proposal distributions that 
approximate the within-model full conditional distributions of the 
model parameters. In this way, their approach mimics aspects of the 
Gibbs sampler. 
However, their approximate full conditional distributions are 
derived from auxiliary within-model MCMC algorithms, 
requiring an extra Markov chain for each rank $K$ to be considered.
As a result, this approach requires  pre-specification of the values 
of $K$,  leaving open 
 the possibility that 
computational effort is spent on values of $K$ with negligible
posterior probability, or that values of $K$ with high posterior 
probability are overlooked
altogether. 

In the analysis of relational data such as social and biological 
networks, 
the row heterogeneity and column heterogeneity of $\m Y$ are 
of equal interest, thus motivating a singular value decomposition 
approach to the data analysis as opposed to a factor analysis. 
The purpose  of this paper is to provide the necessary calculations 
for
model selection, estimation and inference for statistical models based on 
the singular 
value decomposition. 
The results in the following sections provide a means of 
Markov chain Monte Carlo approximation to 
the posterior distribution of the rank and values of the matrix $\m M$. 
This MCMC  algorithm 
is based on  Gibbs sampling
and, unlike  the algorithm of Lopes and West (2002), 
requires no
auxiliary runs
or pre-specification of the values of $K$. 
Additionally, this model-based method allows for
estimation in the presence of missing data or replications, 
and can be incorporated into a generalized linear modeling
framework to allow for the analysis of a variety of data types.
For example, Section 6 discusses a model extension  that allows 
for the analysis and prediction  of
binary relational data such as social or biological networks. 

In the next section we discuss prior distributions for $\{\m U,\m D, \m V\}$ 
given a fixed rank $K$, 
and  show how the uniform distribution for $\m U$ (the invariant measure on 
the Steifel manifold) 
may be 
specified in terms of the full conditional distributions of its column 
vectors. 
Section 3 presents a Gibbs sampling scheme for parameter estimation 
when the rank of $\m M$ is specified. 
In the case of unspecified rank, 
estimation 
can be achieved via 
a prior distribution
which allows the diagonal elements of $\m D$ to each be zero with non-zero 
probability. In Section 4 we consider posterior 
inference under such a prior distribution, 
and 
develop 
a Markov chain Monte Carlo algorithm which moves between 
models with different ranks. This algorithm is constructed via a Gibbs sampling scheme 
which samples each singular value $d_j$ from its conditional 
distribution. This is done marginally over $\m U_{[,j]}$ and $\m V_{[,j]}$, 
and requires a complicated but manageable integration. 
Section 5 presents a small simulation study that 
examines the sampling properties of the Bayesian procedure. It 
is shown that 
the procedure
is able to estimate the true rank of $\m M$ reasonably well for a variety of
matrix sizes, and the squared error of the Bayes estimate $E(\m M|\m Y)$
is typically much lower than that of the least squares estimator.
In Section 6
the singular value decomposition model is extended 
to accommodate non-normal data via a generalized linear modeling 
framework, which is then used in an 
example analysis of  binary relational data. 
A discussion follows in Section 7.

\section{The SVD model and prior distributions}
As described above, our model for an $m\times n$ data matrix 
is  $\m Y=\m M+\m E$, where $\m M$ is a rank $K$ matrix 
 and 
$\m E$ is a matrix of i.i.d.\ mean-zero normally-distributed noise. 
We induce a prior distribution on the matrix $\m M$  by way 
of a prior distribution on the components of its  singular value decomposition
$\{\m U,\m D, \m V\} $. 

For a given rank $K$, we can take $\m U$ to be an $m\times K$ orthonormal 
matrix. 
The set of such matrices is called the Steifel 
manifold and is denoted $\mc V_{K,m}$. 
A natural, non-informative prior distribution for $\m U$ is the uniform 
distribution on $\mc V_{K,m}$, which is the unique probability measure 
on $\mc V_{K,m}$ that is invariant under left and right orthogonal transformations. 
As discussed in Chikuse (2003, Section 2.5), 
a sample $\m U$ from the uniform distribution on the Steifel manifold 
$\mc V_{K,m}$ may be obtained  by first sampling an $m\times K$  
matrix $\m X$ of independent standard normal random variables and then 
setting $\m U = \m X (\m X'\m X)^{-1/2}$. 
Although this construction is straightforward, it doesn't explicitly 
specify conditional distributions of the form 
$p(\m U_{[,j]}|\m U_{[,-j]})$, which are quantities that will be required for 
the estimation procedure outlined in Section 3. 
We now derive these conditional distributions  via a new 
iterative method of generating samples from the uniform distribution 
on $\mc V_{K,m}$. 

Let  $\m U_{[,A]}$ denote the columns of $\m U$ corresponding to a
 subset of column labels $A \subset\{ 1,\ldots, K\}$,
and let
$\m N_A$ be any $m\times (m-|A|)$ matrix whose columns form an
orthonormal basis for the null space of $\m U_{[,A]}$.
A random $\m U \in \mc V_{K,m}$ can be constructed as follows:
\begin{description}
\item{1.} Sample $\m u_1$ uniformly from the unit $m$-sphere and
          set $\m U_{[,1]}=\m u_1$;
\item{2.} Sample $\m u_2$  uniformly from the unit $(m-1)$-sphere and
          set $\m U_{[,2]}= \m N_{ \{1\}} \m u_2$;
\item $\vdots$
\item{$K$.} Sample $\m u_{K}$  uniformly from the unit $(m-K+1)$-sphere and
  set $\m U_{[,K]}= \m N_{\{1,\ldots, K-1\}} \m u_{K}$.
\end{description}
By construction this procedure 
generates an $m\times K$ matrix  $\m U$
having orthonormal columns. The following result also holds:
\begin{prop}
The probability distribution  of $\m U$  is the uniform
probability measure on $\mc V_{K,m}$.
\end{prop}
A proof is provided in the Appendix. 
Since this probability distribution is invariant under
left and right orthogonal transformations of $\m U$
(see, for example, Chikuse 2003), 
it follows that the rows and columns of $\m U$
are exchangeable. 
As a result,
the conditional distribution of $\m U_{[,j]}$ given
any subset $A$ of columns of $\m U$
is equal to the distribution of $\m N_{A} \m u_j$, where
$\m u_j$ is  uniformly distributed on the $(m-|A|)$-sphere.
This fact facilitates the Gibbs sampling of the columns of $\m U$
and $\m V$ from their full conditional distributions, as
described in Section 3. For simplicity we proceed with posterior 
calculations using this uniform prior, even though
the full conditionals derived in Section 3
 indicate that a more general, informative 
conjugate family could be used with no additional computational 
difficulty.

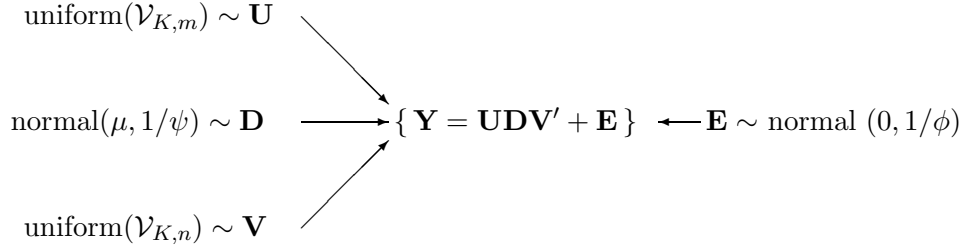
\begin{figure}
\begin{center}
\begin{picture}(230,100)(-30,0)
\put(-80,10){uniform($\mathcal V_{K,n}) \sim \m V$ }
\put(-80,90){uniform($\mathcal V_{K,m}) \sim \m U$ }

\put(-85,50){normal$(\mu ,1/\psi) \sim \m D$ }

\put(25,13){\vector(1,1){33} }
\put(25,53){\vector(1,0){33}} 
\put(25,93){\vector(1,-1){33}}

\put(60,50){$\left \{ \,  \m Y=\m U\m D\m V' +\m E\,  \right  \} $}

\put(178,50){$\m E \sim $ normal $(0,1/\phi) $}

\put(176,53){\vector(-1,0){16}}
\end{picture}
\end{center}
\caption{A graphical representation of the model }
\label{fig:graphmod}
\end{figure}

For a given rank $K$, 
the non-zero singular values $\{ d_1,\ldots, d_K\}$ which make up the diagonal
of $\m D$ determine the magnitude of the mean matrix, in that
$||\m M||^2 = \sum_{k=1}^K d_k^2$.
We model these non-zero 
values as being samples from a normal population with mean
$\mu$ and precision (inverse-variance) $\psi$.
For reasons of conciseness we discuss posterior calculations only for conjugate prior distributions
on these parameters, which include  a
normal distribution  with mean $\mu_0$ and variance $v_0^2$ for
$\mu$, and  a  gamma$(\eta_0/2,\eta_0\tau_0^2/2)$ distribution for
$\psi$, parameterized so that
$\psi$ has expectation $1/\tau_0^2$. 
Other potentially useful prior distributions, along with choices for 
hyperparameters,  are discussed in Section 5. 
This parameterization of the singular values differs slightly from that of the usual singular value decomposition, in that 
the values  $\{d_1,\ldots, d_K\}$  are not restricted to be non-negative here. 
A model enforcing this restriction is possible, but adds a small 
amount of computational difficulty without any modeling benefit (if $\m A$ is 
a diagonal matrix of $\pm 1$'s, then $p(\m Y | \m U,\m D,\m V)  =  
 p(\m Y | \m U\m A,\m A \m D,\m V)  $). 
Finally, the elements of $\m E$ are modeled as i.i.d.\ normal random 
variables with 
mean zero and variance $1/\phi$. 
The prior distribution for the precision $\phi$ 
is taken to be  gamma$(\nu_0/2,\nu_0\sigma_0^2/2)$. 
A graphical representation of the model and  parameters is given in Figure 
 \ref{fig:graphmod}. 
Choices for hyperparameters $\{( \mu_0,v_0^2),( \eta_0,\tau_0^2),(\nu_0 ,\sigma_0^2)\}$ are discussed in Section 5. 

\section{Gibbs sampling for the fixed-rank model}

A Markov chain 
with $p(\m U,\m D,\m V,\phi,\mu,\psi|\m Y,K)$ as its stationary distribution 
can be constructed via a Gibbs sampling procedure, which 
iteratively samples $\phi, \mu,\psi$ and the columns of 
$\m U$, $\m D$ and $\m V$ 
from their 
full conditional distributions. 
These samples can be used to approximate the joint posterior 
distribution and estimate posterior quantities of interest 
(see, for example, Tierney 1994). 

The full conditional distributions for $\phi,\mu ,\psi$  and the elements 
of $\m D$ are standard
and are provided below without derivation. 
Less standard are the full conditional distributions 
of the columns of $\m U$ and $\m V$. 
To derive these, 
consider the form of $p(\m Y|\m U,\m D,\m V, \phi)$ 
as a function of $\m U_{[,j]},\m V_{[,j]}$ and $d_j\equiv \m D_{[j,j]}$. 
Letting 
$ \m E_{-j}  = \m Y -\m U_{[,-j]}\m D_{[-j,-j]}\m V_{[,-j]}'$,  we have
\begin{eqnarray*}
|| \m Y - \m U\m D\m V'||^2 &= & 
 || \m E_{-j}  - d_j \m U_{[,j]} \m V_{[,j]}'||^2 \\ 
&=& || \m E_{-j} ||^2 - 2 d_j \m U_{[,j]}'\m E_{-j} \m V_{[,j]}  + 
  || d_j \m U_{[,j]} \m V_{[,j]}'||^2  \\
 &=& || \m E_{-j} ||^2 - 2 d_j \m U_{[,j]}'\m E_{-j} \m V_{[,j]}  +
   d_j^2. 
\end{eqnarray*}
It follows that  $p(\m Y|\m U,\m D,\m V, \phi)$ can be written 
\begin{equation}
p(\m Y|\m U,\m D,\m V, \phi)= 
  \left ( \frac{\phi}{2 \pi}\right )^{mn/2} 
  \exp \{ -\frac{1}{2} \phi ||\m E_{-j}||^2 +  \phi d_j 
           \m U_{[,j]}'\m E_{-j} \m V_{[,j]}   - \frac{1}{2}\phi d_j^2 \} . 
\label{eq:pyj}
\end{equation}
Recall that conditional on $\m U_{[,-j]}$, $\m U_{[,j]} \stackrel{d}{=}
\m N^u_{ \{-j\} } \m u_j$, where  $\m N^u_{ \{-j\} }$ is a basis for the 
null space of columns of $\m U_{[,-j]}$ and 
$\m u_j$ is uniform on the $m-(K-1)$-sphere. From (\ref{eq:pyj}), 
we see that the  full conditional distribution of $\m u_j$ is proportional 
to $\exp\{\m u_j'(\phi d_j \m N^{u \prime}_{\{-j\}} \m E_{-j} \m V_{[,j]})\}$. 
This is a von Mises-Fisher distribution on the $m-(K-1)$-sphere 
with parameter 
$\phi d_j \m N^{u \prime}_{\{-j\}} \m E_{-j} \m V_{[,j]}$. 
A sample of $\m U_{[,j]}$ from its full conditional distribution 
can therefore be generated by sampling $\m u_j$
from the von Mises-Fisher distribution 
and then setting $\m U_{[,j]} =\m N^u_{\{-j\}} \m u_j$. The full conditional 
distribution of $\m V_{[,j]}$ is derived similarly. 
In general, the von Mises-Fisher distribution on the $p$-sphere 
with  parameter $\bm \mu\in \mathbb R^p$ has density 
$c_p(||\bm \mu||) \exp\{ \m  u'\bm \mu\}$ and is denoted vMF$(\bm \mu)$, and 
the uniform distribution on the sphere is denoted 
vMF$(\bm 0 )$. The normalizing constants for these two cases are
\[ c_p( \kappa )  = (2\pi)^{-p/2}
 \frac{ \kappa ^{p/2-1}}{I_{p/2-1}(\kappa)} \mbox{ for $\kappa>0$,}  \ \ \
c_p(0)  =
 \frac{ \Gamma(p/2) }{2 \pi^{p/2}}  \mbox{ for $\kappa=0$,}
 \]
where $I_{\nu}(x)$ is the modified Bessel function of the first kind
(see Section 9.6 of Abramowitz and Stegun, 1972). 
{\sf R}-code for sampling from this distribution is provided at my
website.

Summarizing these results, a Markov chain with the desired 
stationary distribution can be constructed by iterating the 
following procedure:

\begin{itemize}
\item For $j\in \{1,\ldots, K\}$, 
\begin{itemize}
\item  sample $(\m U_{[,j]}|\m Y,\m U_{[,-j]},\m D , \m V,\phi)  \stackrel{d}{=} \m N^u_{\{-j\}} \m u_j $, where $\m u_j\sim$ vMF$(\phi d_j \m N^{u\prime}_{-j} \m  E_{-j}\m  V_{[,j]} );$
\item  sample $(\m V_{[,j]}|\m Y,\m U, \m D, \m V_{[,-j]},\phi)  \stackrel{d}{=} \m N^v_{\{-j\}} \m v_j $, where $\m v_j\sim$ vMF$(\phi d_j  \m U_{[,j]}'\m  E_{-j}\m N^v_{\{-j\}} );$
\item  sample $(d_j| \m Y,\m U, \m D_{[-j,-j]}, \m V, \phi, \mu ,\psi)\sim $ normal $[ (\m U_{[,j]}'\m  E_{-j} \m V_{[,j]}\phi+ \mu \psi ) /(\phi+\psi)   ,1/(\phi+\psi )]$; 
\end{itemize}
\item sample $(\phi|\m Y,\m U,\m D,\m V)\sim$ gamma$[ (\nu_0 + mn)/2,
(\nu_0 \sigma_0^2 + ||\m Y-\m U\m D\m V'||^2  )/2 ]$;
\item  sample $(\mu |\m D,\psi)\sim$ normal$[(\psi \sum d_j + \mu_0/v_0^2)/(\psi K+ 1/v_0^2), 
     1/(\psi K+ 1/v_0^2  ) ]$
\item  sample $(\psi |\m D,\mu)\sim$ gamma$[ (\eta_0 + K)/2,
(\eta_0 \tau_0^2 +  \sum (d_j-\mu)^2 )/2 ]$;
\end{itemize}

\section{ The variable-rank model}

\subsection{Prior distributions}
In this section we extend the model of Section 2 to the case where
the rank $K$ is to be estimated. This requires comparisons between models 
with parameter spaces of different dimension. 
Two standard ways of viewing such problems
are as follows:
\begin{itemize}
\item Conceptualize a different parameter space for each value of
$K$, i.e., conditional on $K$, the mean matrix is $\m U \m D \m V'$ where the
dimensions of $\m U,\m D$ and $\m V$ are $m\times K$, $K\times K$ and $n \times K$
respectively.
\item Parameterize $\m U,\m D$ and $\m V$ to be of dimensions
 $m\times n$, $n\times n$ and $n \times n$, but allow for columns of
these matrices to be
identically zero. In this parameterization, $K=\sum_{j=1}^n 1(d_j\neq 0)$.
\end{itemize}
Each of these two approaches has its own notational and
conceptual hurdles, and which one to present is to some extent
 a matter of style (see Green 2003 for a discussion). 
Given a prior distribution on $K$, the first approach can be formulated by
using the prior distributions of Section 2 as the conditional distributions
of $\m U,\m D$ and $\m V$ given $K$.
The  second  approach can be made equivalent to
the first as follows: 
\begin{enumerate}
\item Let $\bf \tilde  U,\tilde  D,
\tilde V$ have the prior distributions described 
in Section 2 with $\tilde K=n$; 
\item Let $\{ s_1,\ldots , s_n\}\sim  p(K=\sum s_j )\times
     { n  \choose \sum s_j }^{-1}$, where each $s_j\in\{0,1\}$; 
\item Let $\m S={\rm diag} \{s_1,\ldots, s_n\}$. Set
    $\m U=\bf \tilde U  S,   D=\tilde  D S, V=\tilde  V S$, $K=\sum s_j$. 
\end{enumerate}
Parameterizing a set of nested models with binary variables has 
been a useful technique in a variety of contexts, including 
variable selection in regression models (Mitchell and Beauchamp 1988).
We continue with this formulation because it allows for the 
construction of a relatively straightforward Gibbs sampling scheme 
to generate samples from the  posterior distribution.

The matrices $\m U,\m D$ and $\m V$ described in 1, 2 and 3 above
are exchangeable under simultaneous permutation of their columns.
It follows from Proposition 1 that, 
conditional on $s_1,\ldots, s_n$, 
the non-zero columns of $\m U$ and $\m V$
are random samples from the uniform distributions on
$\mc V_{\sum s_j,m}$ and $\mc V_{\sum s_j,n}$ respectively, and that
conditional on $\{ s_j=1,\m U_{[,-j]},\m V_{[,-j]} \}$,
\[ \m U_{[,j]}\stackrel{d}{=} \m N^u_{\{-j\}} \m u ,  \ \ \ \ \ \
\m V_{[,j]}\stackrel{d}{=} \m N^v_{\{-j\}} \m v, \mbox{    where } \]
\begin{itemize}
\item $\m N^u_{\{-j\}}$ and $\m N^v_{\{-j\}}$ are orthonormal bases for the
null spaces of $\m U_{[,-j]}$ and $\m V_{[,-j]}$;
\item $\m u$ and $\m v$ are uniformly distributed on the
     $(m-\sum s_j+1 )$- and $(n-\sum s_j+1)$-spheres.
\end{itemize}
This property will facilitate posterior sampling of the columns of 
$\m U$, $\m D$ and $\m V$, as described in the next subsection.

\subsection{Posterior estimation}
Let $\bm \Theta = \{ \m U, \m D, \m V\}$, 
    $\bm \Theta_j = \{\m U_{[,j]}, d_j, \m V_{[,j]} \}$ 
  and $\bm \Theta_{-j} = \{  \bm \Theta_k : k\neq j\}$ . 
In this subsection we derive the full conditional distribution
of $\bm \Theta_j$ given $\{ \m Y, \bm \Theta_{-j}, \phi, \mu,\psi\} $
under the model described in the previous subsection. 
The prior 
and full conditional distributions of $\phi,\mu$ and $\psi$ remain 
unchanged from 
Section 2. 
The full conditional distributions can be used in 
a Gibbs sampling scheme to generate approximate samples from 
$p(\m U,\m D, \m V ,\phi, \mu,\psi|\m Y)$.

Under the model and parameterization described above, 
the components of $\bm \Theta_j$ are either all zero or 
have a distribution as described in Section 2. 
To sample $\bm \Theta_j$, we first sample whether or not 
the components are zero, and if not, sample the non-zero values. More 
specifically,  sampling $\bm \Theta_j$ from its full conditional 
distribution can be achieved as follows:
\begin{enumerate}
\item Sample from $( \{ d_j=0\}, \{ d_j \neq 0\}) $ conditional on $\m Y,\bm \Theta_{-j},\phi,\mu,\psi$ . 
\item If $\{ d_j = 0\} $ is true, then set 
  $d_j, \m U_{[,j]}$ and $\m V_{[,j]}$ all equal to zero. 
\item If $\{ d_j \neq 0\}$ is true, 
\begin{enumerate}
\item sample $d_j| \m Y,\bm \Theta_{-j},\phi,\mu,\psi, \{ d_j\neq 0\}$;
\item  sample $\{\m U_{[,j]},\m V_{[,j]}\}|\m  Y,\bm \Theta_{-j},\phi,d_j$. 
\end{enumerate}
\end{enumerate}
The steps 1, 2, and 3 above constitute a draw from
$p( \bm \Theta_{j}| \m Y,   \bm \Theta_{-j}, \phi,\mu,\psi )$.
The first step requires calculation of the odds:
\begin{equation}
  {\rm odds}(d_j\neq 0  |  \m Y,\bm \Theta_{-j},\phi,\mu,\psi ) =  
 \frac{ p(d_j \neq 0 | \bm \Theta_{-j})} 
      { p(d_j=0 | \bm \Theta_{-j} ) }  \times
\frac{ p( \m Y | \bm \Theta_{-j}, d_j\neq 0,\phi, \mu,\psi)  } 
   { p(  \m Y | \bm \Theta_{-j} , d_j=0,\phi,\mu,\psi)}
\label{eq:dodds}
\end{equation}
The first ratio is simply the prior conditional odds of $\{d_j\neq 0\}$ and can be derived 
from the prior distribution on the rank $K$. 
The second term in (\ref{eq:dodds}) can be viewed as a Bayes factor, 
evaluating the evidence in the data for additional structure in $E[\m Y]$
beyond that provided by $\bm \Theta_{-j}$. 
Recall from the previous section that 
$\m Y - \m U \m D \m V' 
 = \m E_{-j} -   d_j\m U_{[,j]} \m V_{[,j]}'$, and so we can write
\begin{eqnarray}
p(\m Y|\m U,\m D,\m V, \phi ,\mu,\psi) &=&
\left [   \left ( \frac{\phi}{2\pi} \right )^{mn/2} \exp\{ -\frac{1}{2}\phi ||\m E_{-j}||^2 \} \right ]
     \exp\{ -\frac{1}{2} \phi d_j^2 \}
     \exp\{ \phi d_j \m U_{[,j]}' \m E_{-j} \m V_{[,j]} \}  \nonumber \\  
 &=& 
 p(\m Y|{\bm \Theta_{-j}}, d_j=0,\phi ,\mu,\psi)  \times 
    \exp\{ -\frac{1}{2} \phi d_j^2 \} \times
     \exp\{ \phi d_j \m U_{[,j]}' \m E_{-j} \m V_{[,j]} \} \label{eq:py} 
\end{eqnarray} 
The 
first term in (\ref{eq:py})  is equal to the denominator of the Bayes factor, 
 and  is simply a product of normal densities
with the elements of $\m Y$ having means given by 
$\m U_{[,-j]} \m D_{[-j,-j]} \m V_{[,-j]}'$ and equal variances $1/\phi$. 
The numerator of the Bayes factor can be obtained by 
integrating  (\ref{eq:py}) over $\bm \Theta_j$ with respect to its conditional 
distribution given  $\mu, \psi, 
\bm \Theta_{-j}$ and $\{d_j\neq 0\}$.  
Integrating first with respect to $\m U_{[,j]},\m V_{[,j]}$, we 
need to calculate 
$ E [ \exp\{ \phi d_j \m U_{[,j]}' \m E_{-j} \m V_{[,j]} \}   |
         \bm \Theta_{-j},d_j ] $. 
Let 
 $\tilde m = m - \sum_{k\neq j} \{ d_k\neq 0\} $ 
and $\tilde n = n-\sum_{k\neq j} \{ d_k\neq 0\}$. 
Recall that conditional on $\bm \Theta_{-j}$, 
$\m U_{[,j]} \stackrel{d}{=} \m N_{\{-j\}}^u\m u$ and $\m V_{[,j]}
\stackrel{d}{=} \m N_{\{-j\}}^v \m v$ 
where
$\m u$ and $\m v$ are uniformly distributed on the $\tilde m$-
and  $\tilde n$-spheres. Letting ${\bf \tilde E} = \m N^{u \prime}_{\{-j\}}
  {\m E}_{-j} {\m N}^v_{\{-j\}} $, 
the required expectation can therefore
be rewritten as 
${\rm E}_{\m u\m v} [ \exp\{ \phi d_j \m u' \bf \tilde E \m v \}  ]$. 
This expectation is non-standard, and is derived in the appendix. 
The result gives: 
\begin{equation}
p(\m Y|\bm \Theta_{-j}, \phi ,\mu,\psi,d_j) =
p( \m Y|\bm \Theta_{-j}, \phi ,\mu,\psi,d_j=0) \times 
\exp\{ -\frac{1}{2} \phi d_j^2 \} 
 \sum_{l=0}^\infty ||{\bf \tilde E}||^{2l} \phi^{2l}  d_j^{2l}  a_l  
\label{eq:marguv}
\end{equation}
where the sequence $\{a_l\}_0^\infty$ can be computed exactly and is 
given in the appendix.

The calculation of $p(\m Y|\bm \Theta_{-j}, \phi ,\mu,\psi,d_j\neq 0)$ 
is completed by integrating  (\ref{eq:marguv}) over $d_j$ 
with respect to $p(d_j| \bm \Theta_{-j}, \phi ,\mu,\psi,d_j\neq 0)$,  the normal density with mean $\mu$ and precision
$\psi$.
This integration simply requires  calculating the 
even moments 
of a normal  distribution, resulting in 
\begin{equation}
p(\m Y|\bm \Theta_{-j}, \phi ,\psi,d_j\neq 0)=
p( \m Y|\bm \Theta_{-j}, \phi ,\psi,d_j=0) \times
 \sum_{l=0}^\infty ||{\bf \tilde E}||^{2l}  a_l b_l
\label{eq:margduv}
\end{equation}
where  the sequence $\{b_l\}_0^\infty $ is given by 
\[ b_l =  \phi^{2l} \left ( \frac{\psi}{\phi+\psi}  \right )^{1/2} 
            \exp\{ -\frac{1}{2} \mu^2 \psi \phi/(\phi+\psi) \}
    E[ \{  \frac{1}{\sqrt{\phi+\psi}}(Z+\frac{\mu\psi}{\phi+\psi}) \}^{2l} ]
\]
where $Z$ is standard normal. The required moments can be calculated 
iteratively, see for example Smith (1995). 
The conditional odds of $\{d_j\neq 0\}$ is therefore 
\[ {\rm odds}(d_j\neq 0  |  \m Y,\bm \Theta_{-j},\phi,\mu,\psi ) = 
 \frac{ p(d_j\neq 0 | \bm \Theta_{-j} )  }
   { p(d_j=0 |\bm \Theta_{-j} )  } \times \sum_{l=0}^\infty ||{\bf \tilde E}||^{2l}   a_l b_l. 
      \]  
In practice, only a finite number of terms can be used to compute the 
above sums. 
The sum in (\ref{eq:marguv}) can be bounded above 
and below by modified Bessel functions, and the error in a finite-sum 
approximation can be bounded. 
This can also provide a guide as to how many terms to include in 
approximating  (\ref{eq:margduv}). 
Details are given in the 
Appendix.

If  $\{d_j\neq 0\}$ is sampled it is still necessary to sample 
 $d_j,\m U_{[,j]}$ and $\m V_{[,j]}$. 
Multiplying equation (\ref{eq:marguv}) by the prior for $d_j|\{d_j\neq 0\}$, 
the required conditional distribution for $d_j$ is proportional to
\[ p(d_j | \m Y , \bm \Theta_{-j},\phi,\mu,\psi, \{d_j\neq 0\} ) \propto
 e^{-\frac{1}{2}(d_j-\mu )^2 \psi} e^{-\frac{1}{2}d_j^2 \phi} \sum_{l=0}^\infty
   ||{\bf \tilde E}||^{2l}  \phi^{2l}  d_j^{2l} a_l\]
which is  an infinite mixture with the following components:
\begin{itemize}
\item mixture weights: $w_l \propto ||{\bf \tilde E}||^{2l}  a_l b_l  $ 
\item mixture densities:
              $f_l(d) \propto d^{2l} \exp\{ -\frac{1}{2} (d-\tilde \mu)^2 \tilde \psi \} $, where $\tilde \mu= \mu \psi/(\phi+\psi)$ and $\tilde \psi = 
    \phi+\psi$
\end{itemize}
The density $f_l(d)$  is nonstandard, but can be sampled from quite efficiently 
using rejection sampling with a scaled and shifted $t$-distribution 
as the  approximating 
density (the tails of a normal distribution are not heavy enough). 

To sample $\m U_{[,j]}$ and $\m V_{[,j]}$ we first sample 
$\m u$ and $\m v$ from their joint distribution and then 
set $\m U_{[,j]} = \m N_{\{-j\}}^u\m u$ and $\m V_{[,j]}
= \m N_{\{-j\}}^v \m v$. 
Equation (\ref{eq:py}) indicates that the joint conditional density of 
$\{ \m u,\m v\}$ is of the form
\begin{equation}
p(\m u,\m v|\m Y,\bm \Theta_{-j},\phi,\mu,\psi, d_j)  = 
c(\m A) \exp\{ \m u'\m A \m v\}  ,
\label{eq:depuv}
\end{equation}
where $\m A= \phi d_j {\bf \tilde E}$ and $c(\m A)^{-1}= 
c_{\tilde m}(0)^{-1}c_{\tilde n}(0)^{-1} \sum_{l=0}^\infty ||\m A||^{2l} a_l$.
This density defines a joint distribution for two dependent unit vectors.
A similar distribution for dependent unit vectors
was discussed  in Jupp and Mardia (1980), except there the vectors 
were of the same length and the matrix $\m A$ was assumed to be orthogonal. 
Some useful facts about the distribution in (\ref{eq:depuv}) include
\begin{itemize}
\item the conditional distribution of $\m u|\m v$ is
vMF$( \m A\m v)$, and that of
$\m v|\m u$ is vMF$(\m u'\m A)$;
\item the marginal distribution of $\m v$ is proportional to 
       $I_{\tilde m/2-1}(||\m A \m v|| ) / || \m A \m v||^{\tilde m/2-1}$; 
\item the joint density has local maxima at $\{ \pm (\hat {\bf u}_k,\hat {\bf v}_k), k=1,\ldots, \tilde n  \}$ 
where $( \hat {\bf u}_k, \hat {\bf v}_k) $ are the
$k$th singular vectors of $\m A$. 
\end{itemize}
I
have implemented  a number of 
rejection and importance
samplers for this distribution, although
making these schemes efficient is still a work in progress.
A relatively fast approximate method 
that seems to work well for a variety of matrices $\m A$
is to first  sample $( \m u,\m v) $ from  the local modes
  $\{ \pm (\hat {\bf u}_k,\hat {\bf v}_k), k=1,\ldots, \tilde n  \}$      
according to the exact relative probabilities
and then 
use this value as a starting point for a 
small number of Gibbs samples, 
alternately sampling  from $p(\m u|\m A,\m v)$ and $p(\m v|\m A , \m u)$.

We now outline a Gibbs sampling algorithm that moves within
and between models of different dimensions, 
using the conditional distributions derived above and 
in Section 3:
\begin{description}
\item[A. Variable dimension sampler:]
 For  each  $j \in \{ 1,\ldots,n \}$, sample $\bm \Theta_j= \{ 
 \m U_{[,j]}, d_j, \m V_{[,j]} \}:$ 
\begin{itemize}
\item  sample $d_j |\m Y, \bm  \Theta_{-j}, \phi, \mu,\psi$  ; 
\item  sample $(\m U_{[,j]}, \m V_{[,j]}) | \m Y, \bm  \Theta_{-j}, 
 \phi,\mu ,\psi, d_j$. 
\end{itemize}
\item[B. Fixed dimension sampler:]
 For each $\{ j : d_j \neq 0 \}$,
\begin{itemize}
\item  sample $\m U_{[,j]} |\m Y, \bm  \Theta_{-j}, \phi, \mu,\psi, d_j, 
\m V_{[,j]}$; 
\item  sample $\m V_{[,j]} |\m Y, \bm  \Theta_{-j}, \phi, \mu,\psi, d_j, 
\m U_{[,j]}$ ;
\item  sample $d_j|\m Y, \bm  \Theta_{-j}, \phi, \mu,\psi,\m U_{[,j]},\m V_{[,j]}$.
\end{itemize}
\item[C. Other terms:] \ \ 
\begin{itemize}
\item sample $\phi|\m Y , \bm \Theta$ ; 
\item sample $\mu| \m D ,\psi$; 
\item sample $\psi|\m D,\mu $. 
\end{itemize}
\end{description}
The distributions required for the steps in $\bf A$ are outlined in this
section, and steps $\bf B$ and $\bf C$ are outlined in the previous
section. 
Technically, interating steps $\bf A$ and $\bf C$ alone will produce 
a Markov chain with the desired stationary ditribution, but 
adding the steps in $\bf B$ can improve the within-model mixing 
of the Markov chain at a relatively low computational cost. 
Also, by conditioning on whether or not $d_j=0$ for each $j$,
the steps in 
 $\bf B$ can  be viewed 
as Gibbs sampling for all $j\in \{ 1,\ldots, n\}$, 
 not just those for which 
$d_j\neq 0$. 
  {\sf R}-code that implements this algorithm is available 
at my website.

\section{Simulation study}
In this section we examine the sampling properties 
of the estimation procedure with a small simulation 
study. Each dataset in this study was simulated from 
the following model:
\begin{itemize}
\item $\m U\sim $ uniform$(\mc V_{5,m})$, 
      $\m V\sim $ uniform$(\mc V_{5,n})$ ; 
\item $\m D  ={\rm diag}\{ d_1,\ldots, d_5\}$, 
      $\{ d_1,\ldots, d_5\}\sim$ i.i.d.\ 
      uniform$( \frac{1}{2}\mu_{mn},\frac{3}{2} \mu_{mn} )$. 
\item $\m Y = \m U \m D \m V' + \m E $, where $\m E$ is an 
     $m\times n$ matrix of standard normal noise. 
\end{itemize}
For each value of $m$ and $n$, the sampling mean of $\{d_1,\ldots, d_5\}$
was taken to be
$\mu_{mn}=\sqrt{n+m+2\sqrt{nm}}$. Such a value should distribute
the singular values $\{ d_1,\ldots, d_5\}$ near the ``cusp'' 
of detectability: As shown in Edelman (1988), the largest 
singular value of an $m \times n$ matrix $\m E$ of standard 
normal noise is approximately $\mu_{mn}$ for large $m$ and $n$. 

Three-hundred datasets
were generated using the model above, one-hundred for each of the three
sample sizes $(m,n) \in \{ (10,10), (100,10),(100,100) \}$.
These were generated in the {\sf R} statistical computing 
environment using the integers 1 through 100 as random seeds
for each of the three sample sizes. Code to generate these datasets 
is available from my website. 
Prior distributions for the  parameters $\{\phi, \mu, \psi\}$ 
were taken as described above with ``prior sample sizes'' of 
$\nu_0=2$ and  $\eta_0=2$. 
This gives 
exponential prior distributions for $\phi$ and $\psi$. 
The values of 
$\sigma_0^2,\mu_0$ and $\tau_0^2$ 
were derived from  ``empirical Bayes''-type estimates 
obtained by averaging over different ranks as follows: 
\begin{enumerate}  
\item For each $k\in\{0,\ldots, n\}$, 
\begin{enumerate}
\item  Let $\hat{\m U} \hat {\m D} \hat {\m V }' $ be the least-squares 
       projection of $\m Y$ onto the set of rank-$k$ matrices; 
\item Let  
$\hat \sigma_k^2 = || \m Y - \hat{\m U} \hat {\m D} \hat {\m V }'||^2/(nm)$
\item Let $\hat \mu_k =  \sum_{j=1}^k \hat d_{j} /k , \ 
           \hat \tau_k^2 =  \sum_{j-1}^k (\hat d_j -\bar {\hat d})^2/k$. 
\end{enumerate}
\item Let $\sigma^2_0 = \frac{1}{n+1} \sum_{j=0}^{n} \hat \sigma^2_j $, 
          $\mu_0 = \frac{1}{n+1} \sum_{j=0}^{n}\hat \mu_j $, 
          $v_0^2 = \frac{1}{n} \sum_{j=0}^n (\hat \mu_j - \bar {\hat \mu})^2$,
          $\tau_0^2 = \frac{1}{n+1} \sum_{j=0}^{n} \hat \tau^2_j $. 
\end{enumerate}
The resulting prior distributions are weakly centered around 
averages of empirical estimates, where the averaging is over ranks
0 through $n$.  Finally, the prior distribution on the rank $K$ of the mean matrix was taken
to be uniform on $\{0,\ldots, n\}$.

\begin{figure}
\centerline{\includegraphics[bb=5 2 570 421,height=4.75in]{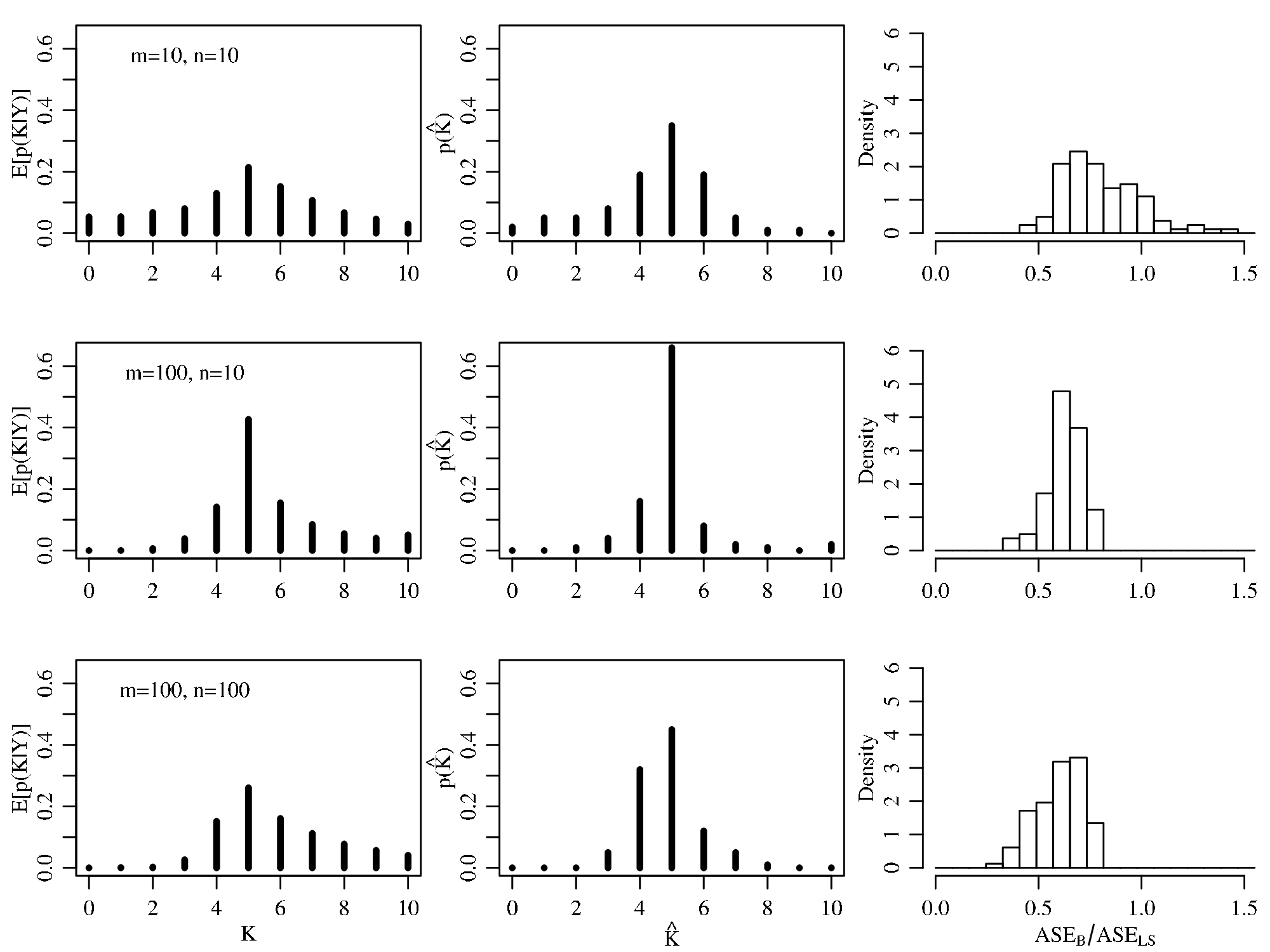}}
\caption{Results of the simulation study. Plots in the first column
give the averages  of $p(K|\m Y)$ over 100 simulated
datasets. The second column gives the empirical distribution
of the posterior mode $\hat K$. The third column gives 
the distribution of the ratio of the squared error of the
Bayes estimate of $\m M$ to that of  the least-squares estimate. }
\label{fig:simres}
\end{figure}  

For each of the $100\times 3$ datasets,  $20,000$ iterations of the
Gibbs sampling scheme described in Section 4.3 were  run  to 
obtain approximate 
samples from the posterior distribution of $\m U \m D\m V'$.  
Parameter values were saved every 10th scan after dropping the 
first 10,000  scans to allow for burn-in, 
resulting in
1000 Monte Carlo samples  per simulation. 
All Markov chains were begun with $K=0$ and 
 $\{ \phi,\mu,\psi \}$ set equal to their prior modes. 
Summaries of the posterior distributions for the three different 
values of $(m,n)$ are displayed in Figure \ref{fig:simres}. 
The first column of each panel plots the MCMC 
approximation to the expected value
of $p(K|\m Y)$ for each value of $(m,n)$. The expectation 
 $E_{Y} [p(K|\m Y) ]$ is approximated by $\frac{1}{100}
 \sum_{s=1}^{100} p(K|\m Y_{(s)})$, 
where $\m Y_{(s)}$ is the $s$th simulated dataset for a given value of $(m,n)$
(for the case $m=n=100$, $p(K|\m Y)$ is plotted only for          
$K\leq 10$, although the distribution extends beyond this value).
These distributions are all peaked around the correct value of $K=5$. 

Also of interest is how frequently the posterior mode
$\hat K = \arg \max p(K|\m Y)$ obtains the true value
of $K$.  This information is displayed  in the second column
of Figure \ref{fig:simres}, which gives the empirical distribution
of $\hat K$ taken over each of the 100 datasets.
As we see, the true value $K=5$ is the most frequent value of the
estimate in each dataset. 
For each dataset we also computed 
three other estimates of $K$:  
$\hat K_l$, the Laplace approximation of Minka (2000); 
$\hat K_c$, the number of eigenvalues of the correlation 
matrix of $\m Y$ that are larger than 1, often suggested  in the 
factor analysis literature;  and 
$\hat K_e$, 
the index of the largest gap in the eigenvalues
of $\m Y'\m Y$, used in machine learning and clustering
(see, for example, Meila and Xu 2003). 
Descriptions of the sampling distributions of these estimators 
are presented in Table {\ref{tab:simstudy}. 
The only case in which the peak of the sampling distribution
for one of these estimators obtained the correct value was 
$\hat K_l$ in the case of $m=100,\ n=10$, although 
the sampling distribution was less peaked around the true 
value than that of the Bayes estimate
($\Pr( \hat K_l =5 | \m Y) = .46$ versus 
  $\Pr( \hat K =5 | \m Y) = .66$). 
 In the case of 
$m=10,n=10$, $\hat K_l$ did poorly, with 
  $\Pr( \hat K_l =1 | \m Y) = .86$. 
Finally, we note that $\hat K_c$ behaved extremely poorly for the 
case $m=n=100$, having a sampling distribution centered around
$K=36$. This is perhaps due to the fact that the asymptotic
behavior of eigenvalues for square matrices are
quite different than that of
rectangular matrices (see Edelman, 1988). 
{\sf R}-code to obtain all estimators for these simulated data 
are available at my website.

\begin{table}
\begin{center}
\begin{tabular}{c||c|c|c|c}
            & \multicolumn{4}{c}{Estimator} \\  \hline
Sample size &  $\hat K $  & $\hat K_e $ & $\hat K_c $ & $\hat K_l $ \\ \hline
 \hline
(10,10)   &  5, .35     & 2, .00     & 4, .14      &  1, .00   \\
(100,10)  &  5, .66    &  3, .09 & 4, .27   &  5, .46    \\
(100,100) &  5, .45       & 6, .21     &  36, .00 &  4, .24   \\
\end{tabular}
\end{center}
\caption{Comparison of model selection procedures.
  Mode of the sampling distribution of $\hat K$ and 
 the probability that $\hat K=5$
for different estimators of the rank of $\m M$. }
\label{tab:simstudy} 
\end{table}

Perhaps of more importance than an estimate of $\hat K$ is an
accurate estimate of $\m M$. In the last column of Figure 
\ref{fig:simres} we compare the error of the model-averaged 
estimate of $\m M$ to that of the least-squares estimate. 
For each simulated dataset the posterior
mean $\hat {\m M} = E[\m M|\m Y]$ was obtained by averaging its
value over the last $10^4$ scans of the Gibbs sampler. The squared
error in estimation, averaged over elements of the mean matrix was calculated as
${\rm ASE}_B = ||\hat {\m M} - \m M||^2/(mn)$ where  
$\m  M$ is the mean matrix that generated the data. This value is
compared to ${\rm ASE}_{LS}$, which is the corresponding average squared error
of the least-squares projection of $\m Y$ onto the space of rank-$\hat K$ matrices. 
The distribution of this ratio is mostly below 1 for
the case $m=n=10$, and strictly below 1 for the other two cases where
there are more parameters to estimate. This corresponds with our intuition:
The model-averaged estimates improve relative to the least-squares
estimates as the number of parameters increases. 
These results indicate that simply obtaining a posterior estimate $\hat K$ of
$K$ and then using the corresponding rank-$\hat K$
least-squares estimate of $\m M$ generally results in an estimate that
can be substantially improved upon by model averaging, at least in
terms of this error criterion. 

As a simple summary of the mixing properties of these Markov chains,
 we examined the convergence and 
autocorrelation of the marginal samples of the error precision $\phi$ using 
Geweke's (1992) $z$-test of stationarity, and by calculating 
the effective sample size of the 1000 Monte Carlo samples saved 
for each chain. 
We declared a Markov Chain simulation a ``success'' if Geweke's 
$z$-statistic did not exceed 2 in absolute value and 
if the effective sample size was at least 100. Based on 
this criterion, the percentage of MCMC simulations that were 
successful was 82\%, 86\% and 91\% for the three sample sizes 
$\{ (10,10),(100,10), (100,100) \}$ respectively. 
The simulation ``failures'' were not examined extensively, 
but we  conjecture that running these chains longer 
is likely to  result improved estimation. 
Of course, in the analysis of a single dataset the usual recommendation is
 to assess the convergence and autocorrelation 
of the Markov chain and to adjust the length accordingly. 

Finally, the performance of the model was examined using two  additional
prior distributions. The posterior analysis was first rerun
using a ``diffuse'' conjugate 
prior distribution,  in which $\phi\sim$ exponential(1), 
 $\psi\sim$ exponential(1) and  $\mu \sim$ normal$(0,1/\psi)$. 
The modes of the sampling distribution of $\hat K$ were
6, 5 and 5 for the sample sizes (10,10), (100,10) and (100,100) 
respectively, with $\Pr(\hat K=5)$ equal to 
.21, .53 and .32 for the three different cases, indicating slightly worse performance 
than the empirical Bayes prior distribution, but better performance 
than the other approaches. A modification to this prior was also 
implemented, in which  $\mu \sim$ 
normal$( \phi^{-1/2}\sqrt{n+m+2\sqrt{nm}}  ,1/\psi)$  and the other 
distributions remained the same. 
This prior distribution focuses the search for non-zero $d_j$'s  to 
values that are as large as the largest singular values of 
normally distributed noise matrices.  Such a prior distribution may be 
appealing in practice, as it requires that factors entering into the model 
have a reasonable magnitude. 
Use of this prior distribution resulted in sampling distributions 
for $\hat K$ having modes of 5 for each of the sample sizes, with 
$\Pr(\hat K=5)$ equal to 
.32, .55, .42. 
A more complicated alternative to these  prior distributions
would be to have the prior distributions for $\{\phi,\mu,\psi\}$
depend on $K$. For example, given
$K=k$, the prior distributions  for $\{\phi,\mu, \psi\}$ could be based on
 $\{\hat \sigma^2_k,\hat \mu_k,\hat \tau^2_k\}$.
Such  prior distributions would require some minor modifications to
the variable-dimension sampler outlined in the previous section.

\section{Extension and example: analysis of binary relational data  }
A potentially useful extension of the model described above is to a class 
of generalized bilinear models  of the form
\begin{eqnarray*} 
\theta_{i,j} &=& \bm \beta'\m x_{i,j} + \m u_i'\m D \m v_j + e_{i,j}  \\
 E[y_{i,j}|\bm \Theta]  &=& g^{-1}(\theta_{i,j}   )
\end{eqnarray*}
where $g$ is the link function. Such  models allow for the analysis 
of a variety of data types: For example, 
binary data can be modeled as $y_{i,j}\sim {\rm binary}(
 \frac{\exp\{\theta_{i,j}\}}{1+\exp\{\theta_{i,j}\}} )$ 
and count data as 
$y_{i,j}\sim {\rm Poisson}(
\exp\{\theta_{i,j}\})$. 
Gabriel (1998) considered maximum likelihood estimation for a variant of 
this model in situations where the dimension of  $\m D$ is fixed, and 
Hoff (2005) considered a symmetric version of this model for the analysis 
of social network data. 
Parameter estimation and dimension selection for the above model can be made by sampling from 
a Markov chain generated by a modified version of 
the algorithm of Section 4.3. 
Given current values of  
$\bm \Theta, \bm \beta, \m U, \m D, \m V$, sample new values as follows:
\begin{enumerate}
\item Let $\tilde{\m  Y} = \bm \Theta - \m X \bm \beta  = 
                         \m U \m D \m V' + \m E$. 
     Update $\m U$, $\m D$,  and $\m V$ from their conditional distribution 
      given $\tilde {\m Y}$ 
     as described in Section 4.3. 
\item Let $\tilde{\m  Y} = \bm \Theta -\m U \m D \m V' = 
         \m X \bm \beta +\m E$. Update $\bm \beta$ from its  conditional 
       distribution given $\tilde {\m Y}$ (a multivariate normal distribution).
\item  Sample ${\bm \Theta}^*  =  \m X \bm \beta  +  \m U \m D \m V' +
         {\m E}^* $, where ${\m E}^*$ is  a matrix of normally 
        distributed noise with zero mean  and precision $\phi$.  
       Replace  $\theta_{i,j}$ by  $ \theta_{i,j}^*$  with probability 
 $\frac{ p(y_{i,j}| \theta_{i,j}^* )}{ p(y_{i,j}|  \theta_{i,j} )}
   \wedge 1 $. 
\end{enumerate}

We illustrate the use of such a model and estimation procedure  with an 
analysis of binary relational data between 46 global service 
firms and 55 cities, obtained from the Globalization and World 
Cities study group ({\tt  http://www.lboro.ac.uk/gawc}). 
For these data, $y_{i,j}=1$ if firm $j$ has an office in city $i$ 
and $y_{i,j}=0$  otherwise. Standard practice is to represent 
within-row and within-column homogeneity with effects that are
additive
on the log-odds scale:
\begin{equation}
  \log {\rm odds} (y_{i,j}=1 ) =  \beta + a_i + b_j ,  \label{eq:rceff}
\end{equation}
and so the effects
$\m a= \{a_1,\ldots, a_m\}$ and $\m b=\{b_1,\ldots, b_n\}$ 
constitute a 
rank-two structure. We look for evidence of higher-order structure 
by considering the model 
\begin{eqnarray} 
\log {\rm odds} (y_{i,j}=1 ) &=&  \beta + \gamma_{i,j}   \label{eq:lom} \\
\gamma_{i,j} &=& {\m u_i}' \m D {\m v_j} + e_{i,j}   \nonumber
\end{eqnarray}
The rank-two structure of model (\ref{eq:rceff})
is easily incorporated into (\ref{eq:lom}) by fixing 
$\m U_{[,1]} = \frac{1}{\sqrt m } \m 1_{m\times 1}$ and 
$\m V_{[,2]} = \frac{1}{\sqrt n } \m 1_{n\times 1}$  and 
modeling $d_1$ and $d_2$ to be non-zero with probability 
1. The additive city and firm effects are then given by 
$\m a= d_2 \m U_{[,2]}$ and $\m b=d_1 \m V_{[,1]}$ 
respectively. Note that any remaining effects represented 
by $\m U \m D \m V'$ will be orthogonal to these additive effects, 
and that the mean of the matrix $\m U \m D \m V'$ is identically zero, 
making it 
unaliased with the intercept $\beta$. 
For the remainder of this analysis, the variable $K$ will refer to 
the number of additional non-zero singular values of 
$\m U \m D\m V'$ beyond the additive row and column effects.

We fix the error variance $1/\phi =1$, as this 
scaling parameter is confounded with the magnitude of $\beta$ and 
$\m U \m D \m V'$. For simplicity we use independent normal $(0,100)$
prior distributions for $\beta$ and the non-zero elements of  
$\m D$, and a uniform prior distribution for $K$. 
A Markov chain of length 25,000 was constructed using the 
algorithm described above, starting with $K=0$. 
Mixing across ranks $K$ was rapid as is shown in 
the first panel of Figure \ref{fig:cities.mcmc}, 
which displays values of $K$ every 100th scan of the Markov chain.
Stationarity of the Markov chain in $K$ was not rejected at level 0.05 based on 
Geweke's $z$-test, and
the effective sample size for estimating the posterior 
distribution of $K$ was 472. 
The 
Monte Carlo estimate of $p(K|\m Y)$, shown in the second panel,
gives a posterior mode of $K=6$ and 
suggests strong evidence for structure in the log-odds beyond
that of the additive row and column effects. 
\begin{figure}
\centerline{\includegraphics[bb=2 2 422 152,height=2.25in]{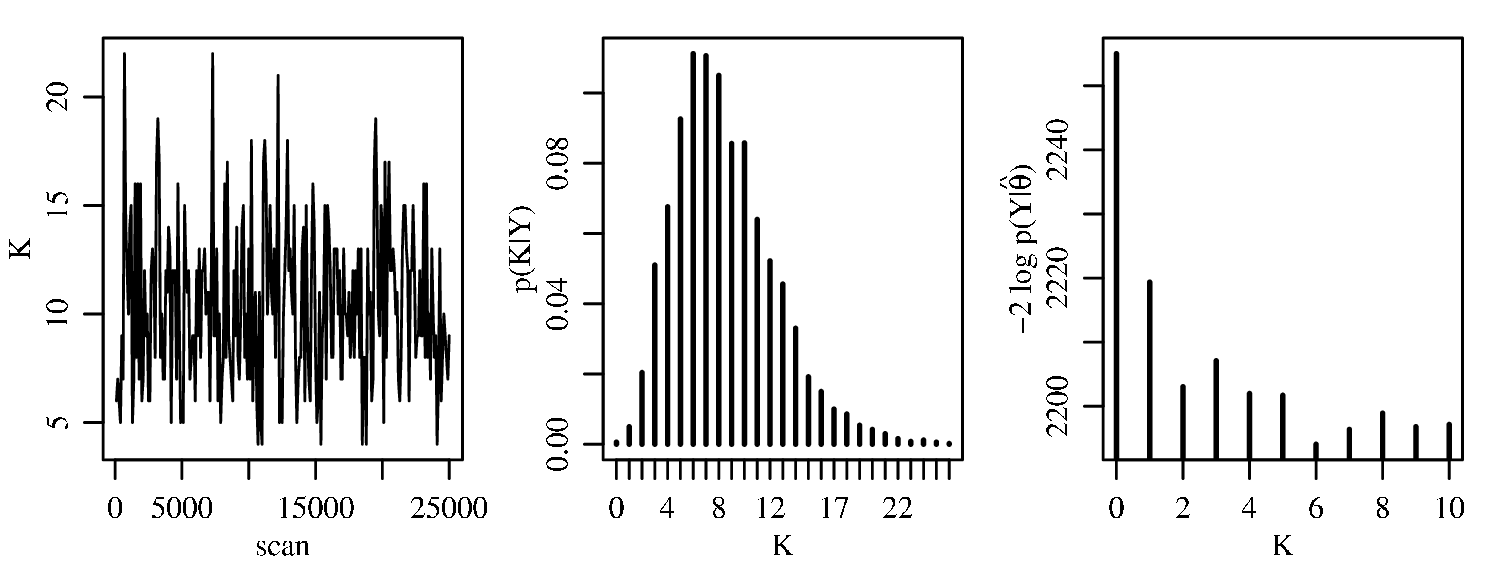}}
\caption{ Posterior estimation of $K$.  The first panel plots values of 
 $K$ every 100th scan of the Markov chain. The second panel plots 
the Monte Carlo estimate of $p(K|\m Y)$. The third panel gives the 
results of a cross-validation evaluation of $K\in \{0,\ldots, 10\}$. }
\label{fig:cities.mcmc}
\end{figure}

One of the practical motivations for selecting an appropriate model dimension
is prediction. 
Many binary social network datasets include missing values, in which 
it is not known whether $y_{i,j}=1$ or $y_{i,j}=0$. 
In such cases it is often  desirable 
to make predictions about missing values based on the observed 
data, and thus to base model selection on predictive performance. 
With this in mind, we compare the above results to the following 
10-fold cross validation procedure:
\begin{enumerate}
\item Randomly split the set of  pairs $\{i,j \}$
       into ten test sets $A_1,\ldots, A_{10}$.
\item For $K=0,1,\ldots, K_{\max} $ :
\begin{enumerate}
\item For $l=1,\ldots, 10$ :
\begin{enumerate}
\item  With the rank fixed at $K$, perform the MCMC algorithm using only
       $\{ y_{i,j}: \{i,j\} \not \in A_{l}\}$,
       but sample values of $\theta_{i,j}$ for all
       ordered pairs.
 \item  Based on the Monte Carlo sample values $\{ \theta_{i,j}^{(1)},\ldots, 
   \theta_{i,j}^{(S)} \}$
       compute the posterior mean  $\hat \mu_{i,j} = 
       \frac{1}{S}\sum_{s=1}^{S} \frac{\exp \{\theta_{i,j}^{(s)}\}}{
       1+ \exp \{\theta_{i,j}^{(s)}\} }$
       for  $\{i,j\}\in A_l $ and the log predictive probability
   lpp($A_l)= \sum_{ \{i,j\} \in A_l } \log p(y_{i,j} | \hat \mu_{i,j} )$.
\end{enumerate}
\item Measure the predictive performance for $K$ as
      LPP$(K) = \sum_{l=1}^{K_{\max}}$ lpp$(A_l)$.
\end{enumerate}
\end{enumerate}
The values of $-2{\rm LPP}(K)$ for $K\in \{0,\ldots,10\}$ are shown in the 
third panel of Figure (\ref{fig:cities.mcmc}).  For the particular 
random partitioning of the data used here, the cross-validation procedure
suggests a model rank of $K=6$, which is the same value as the 
posterior mode of the Bayes solution. 
However, 
a comparison of $N$ values of $K$ using a  ten-fold 
cross validation procedure requires  the construction of 
$10\times N$ separate Markov chains, and further 
requires specification of the values of $K$ to be compared. 
In contrast, the Bayesian procedure
requires only one MCMC run and can  
potentially visit each value of $K\in \{1,\ldots, n\}$. 

Finally we examine some of the patterns in the structure of 
$\m U \m D\m V'$ beyond those of the additive effects. The posterior 
mean of $\m U \m D \m V$, minus the additive effects, was 
obtained by averaging over scans of the Markov chain. 
The first two singular values and 
vectors of this matrix were obtained, and the values of the 
resulting row (city) effects are plotted in 
Figure \ref{fig:cities}.
These values are strongly related to geography: 
U.S.\  cities cluster together, as do cities in Europe,
Latin America and from the Pacific rim.

\begin{figure}
\centerline{\includegraphics[bb=2 0 346 346,height=6in]{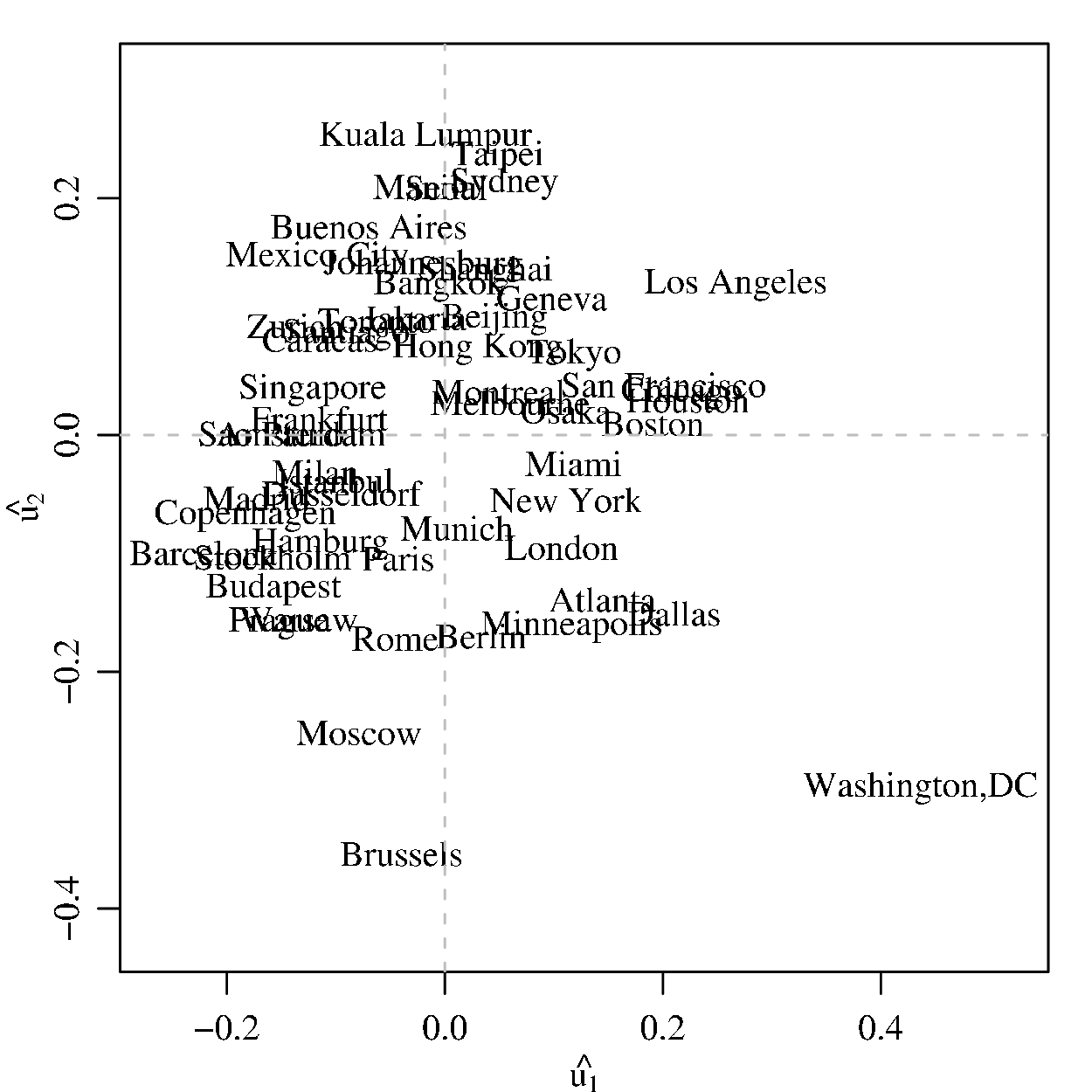}}
\caption{City specific effects: The first two left singular vectors of 
  $ E(\m M|\m Y)$ indicate strong geographic patterns in the data.  }
\label{fig:cities}
\end{figure}

\section{Discussion}
This paper has  presented a model-based
version of the singular value decomposition, thereby extending a 
general data analysis 
tool that has a wide variety of data analysis purposes 
and interpretations. For example, in    
this article it is used for noise reduction and estimation of a mean matrix
(Section 5), as well as prediction and data description (Section 6). 

The approach taken in this paper is to model the data matrix $\m Y$ as 
equal to a reduced-rank mean matrix $\m M$ plus Gaussian 
noise, and to  simultaneously estimate $\m M$  along with its 
rank. 
The approach is Bayesian and  
the estimation procedure, based on Markov chain Monte Carlo, 
allows for a variety of model extensions, such as to the generalized 
bilinear models described in Section 6, 
estimation using replicate 
data matrices and estimation subject to missing data. 
This latter extension may be of particular use in the analysis 
of relational data among a large number of nodes,  where it 
may be too costly to make observations on all possible pairs. 
In such cases, 
the value of $y_{i,j}$ may be missing for many pairs, but
one can  make predictions  
based on estimates $\m u_i,\m D, \m v_j$  obtained from 
the observed data. 
Using this approach  to predict missing 
links in social networks and protein-protein interaction networks is 
one of my current research areas. 
However, for large datasets with 1000 nodes ($10^6$ observations) or more, 
the MCMC scheme in this article becomes prohibitively 
computationally expensive. I am currently studying methods of making 
approximate Bayesian inference for large relational datasets. These 
include Laplace approximations for various components of the MCMC scheme of 
Sections 3 and 4, and using variational methods for approximating 
joint posterior distributions (Jordan et al., 1999). 

Although based on the conceptually straightforward Gibbs sampler, 
complexity 
of the 
full conditional distributions used in the Markov chain of Section 4
 suggest we look
for an alternative procedure. For example,
we could  model only $d_j$  to be zero with non-zero probability,
and sample from its full conditional distribution instead of marginally over
$\m U_{[,j]}$ and $\m V_{[,j]}$.
Unfortunately, an algorithm based on this approach will not mix well across ranks of $\m M$
because
$d_j$, $\m U_{[,j]}$ and  $\m V_{[,j]}$ are dependent to an extreme:
The probability of sampling  $d_j\neq 0$ is essentially
zero unless $\m U_{[,j]}$ and  $\m V_{[,j]}$ are near a pair
of local modes,
but  the
probability of  $\m U_{[,j]}$ and  $\m V_{[,j]}$ being in such a state
is essentially zero if $d_j=0$. 
An alternative approach to sampling from distributions on complicated
sample spaces is to use Metropolis-Hastings type algorithms. 
Based on my initial work on this problem, 
obtaining proposal distributions for these algorithms that achieve even 
minimal acceptance rates 
requires an extreme amount of tuning. 
In contrast, Gibbs sampling for this  model
is possible as shown in this article, 
requires no tuning or pre-specification of model dimensions 
to be considered, 
and, for the examples in this article,
mixes well across  matrices $\m M$ of different ranks.

Computer code and data for all numerical results in this paper 
are available at \\
{\tt www.stat.washington.edu/hoff}.

\appendix

\section{Proof of proposition 1}
We first construct a sample from the uniform distribution on 
$\mc V_{K,m}$ and  then show that it has the desired conditional
distributions. 
Let $\m z_1,\ldots,\m z_K$ be i.i.d.\ multivariate normal
$(0,\m I_{m \times m})$. Let $\m x_1=\m z_1$ and for $j=1,\ldots, K-1$
let
\begin{itemize}
\item $\m X_j = ( \m x_1 \cdots \m x_j )$; 
\item $\m P_j = \m I -  \m X_j(\m X_j'\m X_j)^{-1} \m X_j'$;
\item $\m x_{j+1} = \m P_j \m z_{j+1}$. 
\end{itemize}
Note that $\m P_j$ is the symmetric, idempotent 
projection matrix of $\mathbb R^K$ onto
the null space of $\m X_j$, and so the vectors $\m x_1,\ldots, \m x_{j+1}$ 
are orthogonal. For each $j$, let $\m U_j = \m X_j (\m X_j'\m X_j)^{-1/2}$. 
For $j=K$, we have 
\[ \m X_K'\m X_K = \left ( \begin{array}{c} \m x_1' \\  \m x_2' \\  \vdots
    \\  \m x_K \end{array} \right )
     (  \m x_1 \  \m x_2 \ \cdots \ \m x_K )   =
 \left (  \begin{array}{cccc}|\m x_1|^2 &  0 & \cdots &  0 \\
                            0 & |\m x_2|^2 & \cdots & 0 \\
                            0 & 0 & \cdots & |\m x_K|^2 \end{array} \right  )
\]
and so
\[ \m U_K =   (
\frac{\m x_1}{|\m x_1|},  \
     \frac{\m x_2}{|\m x_2|} , \cdots, \
 \frac{\m x_K}{|\m x_K|}  )
\]
is a matrix of $K$ orthonormal  vectors in $\mathbb R^m$.
The proof will be complete if we can show the following:
\begin{description}
\item[Lemma 1:] The distribution of $\m U_K$ is the uniform distribution on
      $\mc V_{K,m}$.   
\item[Lemma 2:]
$ \m U_{[,k+1]} | \m  U_{k} \stackrel{d}{=}
      \m N_{k} \m u_{k+1}$ where $\m N_k$ is an orthonormal basis for the
      null space of $\m U_k$ and $\m u_{k+1}$ is distributed uniformly on
     the $m-k$ dimensional sphere.
\end{description}

\noindent {\bf Proof of Lemma 1.}
By Theorem 2.2.1 (iii) of Chikuse (2003), an $m\times K$ matrix
of the form $\m X_K (\m X_K'\m X_K)^{-1/2}$ is 
uniformly distributed on $\mc V_{K,m}$ if 
$\m X_K$ is an $m\times K$
random matrix
with rank
 $K$  a.s.\
and having a distribution that is invariant  under left-orthogonal transformations.
We will show left invariance for each $\m X_k$
constructed above by induction.
Let $\m H:\mathbb R^m \rightarrow \mathbb R^m$
be an orthogonal transformation, and note that
$\m H \m X_1 = \m H \m x_1  \stackrel{d}{=} \m x_1 = \m X_1$.
Now suppose $\m H \m X_k \stackrel{d}{=} \m X_k$.
The distribution of $\m H \m X_{k+1}$ is determined by its characteristic function:
\begin{eqnarray*}
E[ \exp\{ i \sum_{j=1}^{k+1} \m t_j' \m H \m x_j \} ] &=&
     E[ \exp\{  i \sum_{j=1}^{k} \m t_j ' \m H\m x_j  \}
     E[ \exp\{ i \m t_{k+1}' \m H \m x_{k+1} \} | \m X_k ] ]  \\
\end{eqnarray*}
Note that $\m t_{k+1}' \m H \m x_{k+1} =
    (\m P_k' \m H' \m t_{k+1})' \m z_{k+1}$, where $\m z_{k+1}$ is a vector
of independent standard normals and independent of $\m X_k$.
Thus the characteristic function can be rewritten as
\begin{equation}
 E[ \exp\{  i \sum_{j=1}^{k}  \m t_j' \m H \m x_j  \}
         \exp\{ -\frac{1}{2} \m t_{k+1}'\m H \m P_k \m P_k'  \m H'\m t_{k+1}\} ] =
E[ \exp\{  i \sum_{j=1}^{k} \m t_j' \tilde {\m x}_j  \}
         \exp\{ -\frac{1}{2} \m t_{k+1}' \tilde  {\m P}_k  \m t_{k+1}\} ]
\label{eq:extk}
\end{equation}
where $\tilde {\m x}_j = \m H \m x_j$ and
\begin{eqnarray*} 
\tilde {\m P}_k =   \m H  \m P_k \m P_k ' \m H'  &=& \m H \m P_k \m H' \\
 &=&      \m H( \m I - \m X_k (\m X_k'\m X_k)^{-1} \m X_k' ) \m H'  \\
 &=&  \m I - \m H \m X_k ( (\m H \m X_k)' (\m H \m X_k) )^{-1}  (\m H \m X_k ) '
\end{eqnarray*}
A similar calculation shows that the distribution of $\m X_j$  is
characterized by
\begin{equation}
 E[ \exp\{ i \sum_{j=1}^{k+1} \m t_j'  \m x_j \} ] =
E[ \exp\{  i \sum_{j=1}^{k} \m t_j'  \m x_j  \}
         \exp\{ -\frac{1}{2} \m t_{k+1}'   \m P_k  \m t_{k+1}\} ], 
 \label{eq:exk}
\end{equation}
By assumption,  $\m X_k \stackrel{d}{=} \m H \m X_k$,
and so $\{ \m x_1,\ldots, \m x_k,\m P_k\} \stackrel{d}{=}
   \{ \tilde {\m x}_1,\ldots, \tilde {\m x}_k, \tilde {\m P}_k \}$ and
the expectations (\ref{eq:extk}) and (\ref{eq:exk}) are equal.
Since the characteristic functions specify the distributions, 
$\m H\m X_{k+1} \stackrel{d}{=} \m X_{k+1}$ and the lemma is proved. 

\noindent {\bf Proof of Lemma 2:}
The vector $\m U_{[,k+1]}$ is constructed as
$\m U_{[,k+1]} = \m P_k \m z_{k+1}/|\m P_k \m z_{k+1}|$.
$\m P_k$ has $m-k$ eigenvalues of one, the rest being zero, giving the
 eigenvalue decomposition
$\m P_k = \m N_k \m N_k'$ where $\m N_k$ is a $m\times (m-k)$ matrix whose
columns form an orthonormal
basis for the null space of $\m U_{k}$.
Substituting in $\m N_k \m N_k'$ for
$\m P_k$ gives
\begin{eqnarray*}
\m U_{[,k+1]} &=&
\frac{ \m N_k \m N_k' \m z_{k+1}  }{|\m N_k \m N_k' \m z_{k+1}|   } \\
& = &\m N_k\frac{\m N_k' \m z_{k+1}   }{ (\m z'\m N_k \m N_k' \m N_k \m N_k' \m z)^{1/2}  } \\
& = &\m N_k\frac{   \m N_k' \m z_{k+1}  }{ (\m z'\m N_k \m N_k'\m z)^{1/2}   } \\
& = & \m N_k\frac{ \m N_k'\m z  }{|\m N_k'\m z|   }
\end{eqnarray*}
Note that
for each $k$, $\m U_k = \m X_k (\m X_k'\m X_k)^{-1/2}$, and so
the projection matrix $\m P_k$ can
be written as
$\m I - \m U_k \m U_k'$, a function of $\m U_k$.
Therefore, given $\m U_k$, $\m U_{[,k+1]}$  is equal in distribution to
$\m N_k$ (a function of $\m U_k$)  multiplied by
$\m N_k'\m z/|\m N_k'\m z|$.
The distribution of $\m N_k'\m z$ can be found via its characteristic function:
For an $m-k$-vector $\m t$
\begin{eqnarray*}
E[ \exp\{ i \m t'(\m N_k'\m z) \} ] &=& E[ \exp\{ i (\m N_k\m t)'\m z \} ] \\
&=& \exp\{ -\frac{1}{2} \m t'\m N_k'\m N_k \m t \} \\
&=& \exp\{ -\frac{1}{2} \m t'\m t \},
\end{eqnarray*}
and so we see that $\m N_k \m z_{k+1}$ is equal in distribution to
an $m-k$-vector of independent standard normal random variables,
and so $\m N_k \m z_{k+1}/|\m N_k \m z_{k+1}|$ is uniformly distributed on
on the $m-k$-sphere.

\section{Expectation of the bilinear form  }
In this section we compute $E[e^{\m u'\m A\m v}]$ for uniformly distributed 
unit vectors $\m u$ and $\m v$ 
and an arbitrary $m\times n$ matrix $\m A$. 
Integrating with respect to $\m v$ can be accomplished by noting that
as a function of $\m v$, $e^{\m u'\m A\m v}$ \
is proportional to the von Mises-Fisher 
distribution on the $n$-sphere $S_n$, 
with parameter $\m u'\m A$:

\begin{eqnarray*}
\int e^{\m u'\m A\m v} p(\m v) \ dS_n(\m v)   &=& 
\int e^{\m u'\m A\m v} c_n(0)  \ dS_n(\m v) \\
&=& \frac{c_n(0)}{c_n(||\m u'\m A||)}\int e^{\m u'\m A\m v} c_n(||\m u'\m A||)  \ dS_n(\m v)  \\
&=& \frac{c_n(0)}{c_n(||\m u'\m A||)}  \\
&=& \Gamma(n/2) (2/||\m u'\m A||)^{n/2-1} I_{n/2-1}(||\m u'\m A||) 
\end{eqnarray*}
where  $I_{\nu} $ is the modified Bessel function of the first kind. 
The series expansion of $I_{n/2-1}(||\m u'\m A||)$ gives 
\[ \Gamma(n/2)(2/||\m u'\m A||)^{n/2-1} I_{n/2-1}(||\m u'\m A||) =
 \sum_{l=0}^\infty  ||\m u'\m A||^{2l} \frac{ \Gamma(n/2) }
 {\Gamma(l+1) \Gamma(l+n/2)4^l }. \]
All the terms in the sum are positive, so 
$ E[e^{\m u'\m A\m v}]$ can be found by replacing 
$ ||\m u'\m A||^{2l}$ with its expectation in the above equation. 
To compute this expectation, let $\m A=\m L\bm \Lambda^{1/2} \m R'$ be the singular value 
decomposition of $\m A$, where $\m L'\m L=\m R'\m R=\m I$ and 
$\bm \Lambda$ is a diagonal matrix of the eigenvalues of $\m A'\m A$. Then 
\begin{eqnarray*}
||\m u'\m A||^2 &=& \m u'\m A\m A'\m u  \\ 
  &=&  \m u'\m L\bm \Lambda^{1/2} \m R' \m R \bm \Lambda^{1/2}\m  L'\m  u \\ 
&=& \m u'\m L\bm \Lambda \m L' \m u \\
&\equiv&  \tilde {\m u}'\bm  \Lambda \tilde{\m u} \\
&=&  \sum_{j=1}^n \tilde u_j^2 \lambda_j,   \\
\end{eqnarray*}
where $\tilde {\m u} = \m L'\m u$. 
We will now identify the distribution of the vector
$\{ \tilde u_1^2,\ldots, \tilde u_n^2\}$. 
Let $\m B= \{\m  L , \m L^\perp\}$ be an orthonormal basis for $\mathbb R^m$. 
Since the uniform distribution on the sphere is rotationally 
invariant, $\m B'\m u$ is equal in distribution to $\m u$, and so 
$\m L'\m u$ is equal in distribution to the first $n$ coordinates of 
$\m u$. Recall that a uniformly distributed vector $\m u$ can be generated by 
sampling $z_1,\ldots, z_m$ independently from a standard normal distribution 
and then dividing each term by $|\sum z_i^2|^{1/2}$. Therefore, 
\begin{eqnarray*}
 \{ \tilde u_1^2,\ldots, \tilde u_n^2\} & \stackrel{d}{=} &
   \frac{ \{ z_1^2, \ldots, z_n^2 \} } { \sum_{j=1}^{m} z_j^2 }  \\
 &=& \left ( \frac{  \sum_{j=1}^{n} z_j^2 } { \sum_{j=1}^{m} z_j^2}  \right )
\left(\frac{ \{ z_1^2, \ldots, z_n^2 \} } { \sum_{j=1}^{n} z_j^2 }  \right )\\
&\stackrel{d}{=}& \theta \m q
\end{eqnarray*}
where 
$\theta\sim $ beta$(n/2,(m-n)/2)$, 
$\m q\sim $ $\rm Dirichlet_n(1/2,\ldots,1/2)$
and $\theta$ and $\m q$ are independent. 
Therefore, $||\m u'\m A||^2 \stackrel{d}{=} \theta \bm \lambda'\m q$, 
where $\bm \lambda$ is the diagonal of $\bm \Lambda$ and 
are the eigenvalues of $\m A'\m A$. The required 
expectation is then 
\[ E[ ||\m u'\m A||^{2l} ] = 
E [  \theta^l ] E [ (\bm \lambda ' \m q )^{l} ] 
\]
The first expectation is given by 
$E [  \theta^l ] =  [\Gamma(n/2+l) \Gamma(m/2) ]/[\Gamma(m/2+l)\Gamma(n/2) ]$.
The second expectation is the $l$th-moment of a Dirichlet average, 
which results in a type of multiple hypergeometric function 
denoted as $R_l(\bm \lambda,\frac{1}{2}\bm 1)$. 
This expectation and its generalizations have been studied by 
Carlson (1977, Chapter 5), Dickey (1983) and others. 
An algorithm for recursively computing $R_1,\ldots,R_l$ exactly 
from a generating function is provided 
in the next section. 

To make the result of the calculation a little more intuitive
let $\tilde {\bm \lambda} = \bm \lambda/\sum \lambda_j \equiv \bm \lambda/||\m A||^2$, 
and make use of the fact that 
 $E[ (\bm \lambda'\m q)^l] =  ||\m A||^{2l}
  E [(\tilde {\bm \lambda}'\m q)^l] $, 
 Combining the
results gives 
\[  E[e^{\m u'\m A\m v} ] = 
   \sum_{l=0}^\infty 
   ||\m A||^{2l}   E[(\tilde {\bm \lambda} '\m q)^l ]
\frac{\Gamma(m/2)  }{\Gamma(m/2+l)\Gamma(1+l) 4^l }  
\equiv 
 \sum_{l=0}^\infty 
   ||\m A||^{2l} a_l  , 
\]
and so we see how the expectation is related to the 
norm of $\m A$  via $||\m A||^{2l}$ and the variability in relative sizes of 
the squared singular values via $E[(\tilde {\bm \lambda}'\m q)^l]$. 

To get bounds on a finite-sum approximation to $E[e^{\m u'\m A\m v}]$, 
note that 
$    \lambda_{\min}^l < E[ ( {\bm \lambda}'\m q)^l ]< 
      \lambda_{\max}^l $
so 
\[  
  \sum_{l=r+1}^\infty  \lambda_{\min}^l  \frac{\Gamma(m/2)  }{\Gamma(m/2+l)\Gamma(1+l) 4^l }    <
  \sum_{l=r+1}^\infty E[ ( {\bm \lambda}'\m q)^l ]   \frac{\Gamma(m/2)  }{\Gamma(m/2+l)\Gamma(1+l) 4^l } <
  \sum_{l=r+1}^\infty  \lambda_{\max}^l  \frac{\Gamma(m/2)  }{\Gamma(m/2+l)\Gamma(1+l) 4^l } \]
The outer sums can be computed as 
\[ \sum_{l=r+1}^\infty  \lambda^l  \frac{\Gamma(m/2)  }{\Gamma(m/2+l)\Gamma(1+l) 4^l }  = 
   \left ( \frac{2}{\sqrt{\lambda}} \right)^{m/2-1} I_{m/2-1}(\sqrt{\lambda}) 
    \Gamma(m/2)    - 
 \sum_{l=1}^r  \lambda^l  \frac{\Gamma(m/2)  }{\Gamma(m/2+l)\Gamma(1+l) 4^l, } 
\]
and so bounds  on 
 $ E[e^{\m u'\m A\m v}] - \sum_{l=0}^{r} ||\m A||^{2l} a_l $
 can be obtained.

\section{Computing the multiple hypergeometric function} 
Let $\m q\sim \ {\rm Dirichlet}_n(\alpha_1,\ldots,\alpha_n)$. 
Carlson (1977, Section 6.6) shows that 
\[ \prod_{i=1}^n (1-t \lambda_i)^{-\alpha_i} =  
   \sum_{l=0}^\infty 
   \frac{ \Gamma(\bm \alpha'\m 1 +l) }{\Gamma(\bm \alpha'\m 1) \Gamma(l+1)}
    t^l    E[ (\bm \lambda'\m q )^l ].  \]
Let $c_l = \frac{ \Gamma(\bm \alpha'\m 1 +l) }{\Gamma(\bm \alpha'\m 1)
   \Gamma(l+1)}
      E[ (\bm \lambda'\m q )^l ]$.  We now show how to calculate $c_{k+1}$ 
based on $c_1,\ldots, c_{k}$. 
Let $f(t) =  \sum_{l=0}^\infty c_l t^l$ be the right-hand side of the equation
and $g(t) =   - \sum_{i=1}^n \alpha_i \log(1-t \lambda_i)$ be the log 
of the left-hand side. 
Taking derivatives with respect to $t$ and evaluating at zero we 
have 
\[ f^{(l)}(0) =  \Gamma(l+1) c_l  ,   \ \ \ \ \ \ 
   g^{(l)}(0) =  \Gamma(l) \sum_{i=1}^{n} \alpha_i  \lambda_i^l . \]
Since $f(t)=e^{g(t)}$, we have 
\[  f^{(k+1)}(0) = \sum_{l=0}^{k} { k \choose l } f^{(l)}(0) g^{(k+1-l)}(0). \]
Plugging the values of $f^{(l)}(0)$ into  the sum gives
\[ c_{k+1} = \sum_{l=0}^{k} \left [ c_l { k\choose l} \frac{ \Gamma(l+1)\Gamma(k+1-l) }
 {\Gamma(k+2)}   \left ( \sum_{i=1}^n \alpha_i \lambda_i^{k+1-l} \right ) 
 \right ] . \]
Simplifying gives 
\[ E[ (\bm \lambda'\m q)^{k+1} ] = 
\sum_{l=0}^{k} \left [ E[(\bm \lambda'\m q)^{l}]\frac{ \Gamma(\bm 1'\bm \alpha +l ) \Gamma(k+1)}
{  \Gamma (\bm 1'\bm \alpha + k+1) }  \left ( \sum_{i=1}^n \alpha_i \lambda_i^{k+1-l} \right ) \right ].  \]
{\sf C}-code with an {\sf R}-interface 
to calculate $\{  E[ (\bm \lambda'\m q)^l : l=0,\ldots, k\}$ 
is available at my website.

\nocite{*}
\bibliographystyle{plain}
\bibliography{bsvd}

\end{document}